\documentclass[12pt,reqno]{amsart}

\title{Minimization of entropy functionals}
\author{Christian L\'eonard}
\date{October 07}

\usepackage{amssymb, amsmath, amsfonts, latexsym, enumerate}

\newtheorem{theorem}[equation]{Theorem}

\newtheorem{proposition}[equation]{Proposition}

\newtheorem{definitions}[equation]{Definitions}

\theoremstyle{remark}
\newtheorem{remark}[equation]{Remark}
\newtheorem{remarks}[equation]{Remarks}

\numberwithin{equation}{section}

\newcommand{\R}{\mathbb{R}}

\newcommand{\1}{\textbf{1}}

\newcommand{\eqdef}{\stackrel{\vartriangle}{=}}
\newcommand{\dom}{\mathrm{dom\,}}
\newcommand{\icordom}{\mathrm{icordom\,}}
\newcommand{\icor}{\mathrm{icor\,}}

\newcommand{\cl}{\mathrm{cl\,}}

\newcommand{\inter}{\mathrm{int\,}}

\newcommand{\lsc}{lower semicontinuous}

\newcommand\seq[2]{(#1_#2)_{#2\ge1}}
\newcommand\Lim[1]{\lim_{#1\rightarrow\infty}}
\newcommand\Liminf[1]{\liminf_{#1\rightarrow\infty}}
\newcommand\Limsup[1]{\limsup_{#1\rightarrow\infty}}
\newcommand{\Fo}{\Phi} %{\Phi_o}
\newcommand{\Fos}{\Phi^*} %{\Phi_o^*}
\newcommand{\Fi}{\Phi}
\newcommand{\Fs}{\Phi^*}
\newcommand{\Fb}{\bar\Phi}

\newcommand{\Ib}{\bar I}

\newcommand{\Il}{I_\lambda}
\newcommand{\Ils}{I_{\lambda^*}}

\newcommand{\Lao}{\Lambda_o}
\newcommand{\Los}{\Lambda_o^*}
\newcommand{\La}{\Lambda}
\newcommand{\Ls}{\Lambda^*}
\newcommand{\Lb}{\overline{\Lambda}}
\newcommand{\To}{T_o}

\newcommand{\UU}{\mathcal{U}}
\newcommand{\LL}{\mathcal{L}}
\newcommand{\YY}{\mathcal{Y}}
\newcommand{\XX}{\mathcal{X}}
\newcommand{\ZZ}{\mathcal{Z}}
\newcommand{\CC}{\mathcal{C}}

\newcommand{\UUo}{\mathcal{U}_o}
\newcommand{\LLo}{\mathcal{L}_o}
\newcommand{\YYo}{\YY_o}
\newcommand{\XXo}{\XX_o}

\newcommand{\YYt}{\widetilde{\YY}}
\newcommand{\YYb}{\overline{\YY}}

\newcommand{\MZ}{M_\ZZ}
\newcommand{\PZ}{P_\ZZ}

\newcommand{\lmax}{{\lambda_{\diamond}}}
\newcommand{\lmaxs}{{\lambda_{\diamond}^*}}
\newcommand{\El}{E_\lmax}
\newcommand{\EEl}{\mathcal{E}_\lmax}
\newcommand{\LLl}{\mathcal{L}_\lmax}
\newcommand{\Ll}{L_\lmax}
\newcommand{\Lls}{L_\lmaxs}
\newcommand{\LlsR}{L_\lmaxs R}

\newcommand{\IZ}{\int_{\ZZ}}

\newcommand{\HF}{$(H_\Phi)$}
\newcommand{\HFi}{$(H_{\Fi 1})$}
\newcommand{\HFii}{$(H_{\Fi 2})$}
\newcommand{\HFiii}{$(H_{\Fi 3})$}
\newcommand{\HT}{$(H_T)$}
\newcommand{\HTi}{$(H_{T1})$}
\newcommand{\HTii}{$(H_{T2})$}
\newcommand{\HC}{$(H_{C})$}

\newcommand{\pc}{\mathrm{P}_C}
\newcommand{\pco}{\mathrm{P}_{C_o}}
\newcommand{\pbco}{\overline{\mathrm{P}}_{C_o}}
\newcommand{\pbc}{\overline{\mathrm{P}}_C}

\newcommand{\dc}{\mathrm{D}_C}
\newcommand{\dbc}{\overline{\mathrm{D}}_C}

\newcommand{\dtc}{\widetilde{\mathrm{D}}_C}

\newcommand{\Po}{\mathcal{P}_o}

\newcommand{\Db}{\overline{\mathcal{D}}}

\newcommand{\NF}{|\cdot|_\Phi}
\newcommand{\NL}{|\cdot|_\Lambda}

\newcommand{\sLE}{\sigma(\LlsR,\El)}

\newcommand{\sUL}{\sigma(\UUo,\LLo)}

\newcommand{\sLUi}{\sigma(\LL,\UU)}

\newcommand{\Co}{C_o}

\newcommand{\lh}{\hat\ell}
\newcommand{\ob}{\bar\omega}

\newcommand{\xh}{\hat x}
\newcommand{\ot}{\tilde \omega}

\newcommand{\Qh}{\widehat{Q}}

\newcommand{\ul}{\langle u,\ell\rangle}

\newcommand{\xo}{\langle x,\omega\rangle}
\newcommand{\yx}{\langle y,x\rangle}
\newcommand{\yt}{\langle y,\theta\rangle}
\newcommand{\ott}{\langle \ot,\theta\rangle}

\begin{document}

% ************** page de garde ******************************

 \address{Modal-X, Universit\'e Paris 10\quad \&  B\^at. G,
 200 av. de la R\'epublique. 92001 Nanterre Cedex, France}
 \address{CMAP, \'Ecole Polytechnique. 91128 Palaiseau Cedex, France}
 \email{christian.leonard@polytechnique.fr}

\begin{abstract}
Entropy functionals (i.e.\! convex integral functionals) and
extensions of these functionals are minimized on convex sets. This
paper is aimed at reducing as much as possible the assumptions on
the constraint set. Dual equalities and characterizations of the
minimizers are obtained with weak constraint qualifications.
\end{abstract}

\maketitle \tableofcontents

% ************* corps du texte ****************************

\section{Introduction}

\subsection{The entropy minimization problem}\label{sec:entropmin}

Let $R$ be a positive measure on a space $\ZZ.$ Take a
$[0,\infty]$-valued measurable function $\gamma^*$ on $\ZZ\times
\mathbb{R}$ such that $\gamma^*(z,\cdot):=\gamma_z^*$ is convex
and \lsc\ for all $z\in \ZZ.$ Denote $\MZ$ the space of all signed
measures $Q$ on $\ZZ.$ The entropy functional to be considered is
defined by
\begin{equation}\label{eq-16}
    I(Q)=
    \left\{
    \begin{array}{ll}
    \IZ \gamma^*_z(\frac{dQ}{dR}(z))\,R(dz) & \mathrm{if \ }
    Q\prec R\\
    +\infty & \mathrm{otherwise}
    \end{array}
\right.,\quad Q\in\MZ
\end{equation}
where $Q\prec R$ means that $Q$ is absolutely continuous with
respect to $R.$ Assume that for each $z$ there exists a unique
$m(z)$ which minimizes $\gamma^*_z$ with
\begin{equation}\label{eq-04}
    \gamma^*_z(m(z))=0,\ \forall z\in\ZZ.
\end{equation}
Then, $I$ is $[0,\infty]$-valued, its unique minimizer is $mR$ and
$I(mR)=0.$
\\
This paper is concerned with the minimization problem
\begin{equation}\label{eq-01}
      \textsl{minimize } I(Q)
    \textsl{ subject to } \To Q\in C,\quad Q\in \MZ
\end{equation}
where $\To: \MZ\to\XXo$ is a linear operator which takes its
values in a vector space $\XXo$ and $C$ is a convex subset of
$\XXo.$

\subsection{Presentation of the results}
\emph{Our aim is to reduce as much as possible the restrictions on
the convex set $C.$} Denoting the minimizer $\Qh$ of
(\ref{eq-01}), the geometric picture is that some level set of $I$
is tangent at $\Qh$ to the constraint set $\To^{-1}C.$ Since these
sets are convex, they are separated by some affine hyperplane and
the analytic description of this separation yields the
characterization of $\Qh.$ Of course Hahn-Banach theorem is the
key. Standard approaches require $C$ to be open with respect to
some given topology in order to be allowed to apply it. In the
present paper, one chooses to use a topological structure which is
designed  for the level sets of $I$ to ``look like'' open sets, so
that Hahn-Banach theorem can be applied without assuming to much
on $C.$

This strategy is implemented in \cite{p-Leo07a} in an abstract
setting suitable for several  applications. It is a refinement of
the standard saddle-point method \cite{Roc74} where convex
conjugates play an important role. The proofs of the present
article are applications of the general results of
\cite{p-Leo07a}.

Clearly, for the problem (\ref{eq-01}) to be attained, $\To^{-1}C$
must share a supporting hyperplane with some level set of $I.$
This is the reason why it is assumed to be closed with respect to
the above mentioned topological structure. This will be the only
restriction to be kept together with the interior specification
(\ref{eq-25}) below.

Dual equalities and primal attainment are obtained under the
weakest possible assumption:
$$
C\cap\To\dom I\not=\emptyset
$$
where $\dom I:=\{Q\in\MZ;I(Q)<\infty\}$ is the effective domain of
$I$ and $\To\dom I$ is its image by $\To.$ The main result of this
article is the characterization of the minimizers of (\ref{eq-01})
in the interior case which is specified by
\begin{equation}\label{eq-25}
    C\cap\icor(\To\dom I)\not=\emptyset
\end{equation}
 where $\icor(\To\dom I)$ is the intrinsic core of $\To\dom I.$ The notion of intrinsic core does not rely on any
topology; it gives the largest possible interior set. For
comparison, a usual form of constraint qualification required for
the representation of the minimizers of (\ref{eq-01}) is
\begin{equation}\label{eq-26}
    \inter(C)\cap\To\dom I\not=\emptyset
\end{equation}
where $\inter (C)$ is the interior of $C$ with respect to a
topology which is not directly connected to the ``geometry'' of
$I.$ In particular, $\inter (C)$ must be nonempty; this is an
important restriction. The constraint qualification (\ref{eq-25})
is weaker.

An extension of Problem (\ref{eq-01}) is also investigated. One
considers an extension $\Ib$ of the entropy $I$ to a vector space
$L_\ZZ$ which contains $\MZ$ and may also contain singular linear
forms which are not $\sigma$-additive. The extended problem is
\begin{equation}\label{eq-27}
     \textsl{minimize } \Ib(\ell) \textsl{ subject to }\To\ell\in C, \quad\ell\in L_\ZZ
\end{equation}
Even if $I$ is strictly convex, $\Ib$ isn't strictly convex in
general so that (\ref{eq-27}) may admit several minimizers. There
are situations where  (\ref{eq-01}) is not attained in $\MZ$ while
(\ref{eq-27}) is attained in $L_\ZZ.$ Other relations between
these minimization problems are investigated by the author in
\cite{p-Leo07d} with probabilistic questions in mind.

\subsection{Literature about entropy minimization}\label{sec-lit-ent}
Entropy minimization problems appear in many areas of applied
mathematics and sciences. The literature about the minimization of
entropy functionals under convex constraints is considerable: many
papers are concerned with an engineering approach, working on the
implementation of numerical procedures in specific situations. In
fact, entropy minimization is a popular method to solve ill-posed
inverse problems.
\\
Rigorous general results on this topic are quite recent. Let us
cite, among others, the main contribution of Borwein and Lewis:
\cite{BL91a}, \cite{BL91b}, \cite{BL91c}, \cite{BL92},
\cite{BL93}, \cite{BLN} together with the paper \cite{TV} by
Teboulle and Vajda. In these papers, topological constraint
qualifications of the type of (\ref{eq-26}) are required. Such
restrictions are removed here.
\\
With a geometric point of view, Csisz\'ar \cite{Csi75, Csi84}
provides a complete treatment of (\ref{eq-01}) with the relative
entropy (see Section \ref{sec-expl-ent}) under the weak assumption
(\ref{eq-25}). The behavior of minimizing sequences of general
entropy functionals is studied in \cite{Csi95}.
\\
By means of a method different from the saddle-point approach, the
author  has already studied in \cite{Leo01b, Leo03}  entropy
minimization problems under affine constraints (corresponding to
$C$ reduced to a single point) and  more restrictive assumptions
on $\gamma^*.$
\\
The present article extends these results.

\subsection*{Outline of the paper}
The minimization problems (\ref{eq-01}) and (\ref{eq-27}) are
described in details at Section \ref{sec-minpb}.  In Section
\ref{sec-prelim}, the main results of \cite{p-Leo07a} about the
extended saddle-point method are recalled. Section \ref{sec-pbc}
is devoted to the extended problem (\ref{eq-27}) and Section
\ref{sec-pc} to (\ref{eq-01}). One presents important examples of
entropies and constraints at Section \ref{sec-expl-ent}.

\subsection*{Notation}
Let $X$ and $Y$ be topological vector spaces. The algebraic dual
space of $X$ is $X^{\ast},$ the topological dual space of $X$ is
$X'.$ The topology of $X$ weakened by $Y$ is $\sigma(X,Y)$ and one
writes $\langle X,Y\rangle$ to specify that $X$ and $Y$ are in
separating duality.
\\
Let $f: X\rightarrow [-\infty,+\infty]$ be an extended numerical
function. Its convex conjugate with respect to $\langle
X,Y\rangle$ is $f^*(y)=\sup_{x\in X}\{\langle x,y\rangle
-f(x)\}\in [-\infty,+\infty],$ $y\in Y.$ Its subdifferential at
$x$ with respect to $\langle X,Y\rangle$ is $\partial_Y
f(x)=\{y\in Y; f(x+\xi)\geq f(x)+\langle y,\xi\rangle, \forall
 \xi\in X\}.$ If no confusion occurs, one writes
$\partial f(x).$
\\
The intrinsic core of a subset $A$ of a vector space is
    $
\icor A=\{x\in A; \forall x'\in\mathrm{aff\,}A, \exists t>0,
[x,x+t(x'-x)[\subset A\}
    $
where $\mathrm{aff\,}A$ is the affine space spanned by $A.$
$\icordom f$ is the intrisic core of the effective domain of $f:$
$\dom f=\{x\in X; f(x)<\infty\}.$
\\
The indicator of a subset $A$ of $X$ is defined by
    \begin{equation*}
    \iota_A(x)=\left\{%
\begin{array}{ll}
    0, & \hbox{if }x\in A \\
    +\infty, & \hbox{otherwise} \\
\end{array}%
\right.,\quad x\in X.
    \end{equation*}
The support function of $A\subset X$ is $\iota_A^*(y)=\sup_{x\in
A}\langle x,y\rangle,$ $y\in Y.$
\\
One writes $I_\varphi(u):=\IZ
\varphi(z,u(z))\,R(dz)=\IZ\varphi(u)\,dR$ and $I=I_{\gamma^*}$ for
short, instead of (\ref{eq-16}).

\section{Presentation of the minimization problems ($\pc$) and ($\pbc$)}\label{sec-minpb}

The problem (\ref{eq-01}) and its extension (\ref{eq-27}) are
introduced. Their correct mathematical statements necessitate the
notion of Orlicz spaces. The definitions of good and bad
constraints are given and the main assumptions are collected at
the end of this section.

\subsection{Orlicz spaces} \label{sec:Orlicz}
\newcommand{\Mr}{E_\rho}
\newcommand{\Lr}{L_\rho}
\newcommand{\Lrs}{L_{\rho^*}}
\newcommand{\Lrsi}{L_{\rho}^s}
\newcommand{\la}{\ell^a}
\newcommand{\lsing}{\ell^s}

To state the minimization problem (\ref{eq-01}) and its extension
correctly, one will need to talk in terms of Orlicz spaces related
to the function $\gamma^*.$
\\
 Let us recall some basic definitions and results.
A set $\ZZ$ is furnished with a $\sigma$-finite nonnegative
measure $R$ on a $\sigma$-field which is assumed to be
$R$-complete. A function $\rho:\ZZ\times\R$ is said to be a
\emph{Young function} if for $R$-almost every $z,$ $\rho(z,\cdot)$
is a convex even $[0,\infty]$-valued function on $\R$ such that
$\rho(z,0)=0$ and there exists a measurable function $z\mapsto
s_z>0$ such that $0<\rho(z,s_z)<\infty .$
\\
In the sequel, every numerical function on $\ZZ$ is supposed to be
measurable.

\begin{definitions}[The Orlicz spaces $\mathcal{L}_\rho, \mathcal{E}_\rho,$ $L_\rho$ and $E_\rho$]
The \emph{Orlicz space} associated with $ \rho $ is defined by
$\mathcal{L}_\rho(\ZZ,R)=\{u:\ZZ\rightarrow\mathbb{R}; \Vert
u\Vert_ \rho < + \infty \}$ where the Luxemburg norm $ \Vert \cdot
\Vert_ \rho $ is defined by $\Vert u\Vert_ \rho= \inf \left\{
\beta
>0\ ;\ \IZ\rho(z,u(z)/\beta)\, R(dz)\leq  1 \right\}.$ Hence,
 $$
 \mathcal{L}_\rho(\ZZ,R)=\left\{u:  \ZZ \rightarrow  \R \ ; \exists \alpha_o>0, \IZ
 \rho \Big(z,\alpha_o u(z)\Big)\,R(dz)<\infty \right\}.
 $$
 A subspace of interest is
 \begin{equation*}
     \mathcal{E}_\rho(\ZZ,R)=\left\{u:  \ZZ \rightarrow  \R \ ; \forall \alpha >0, \IZ
 \rho \Big(z,\alpha u(z)\Big)\,R(dz)<\infty \right\}.
\end{equation*}
Now, let us identify the $R$-a.e.\! equal functions. The
corresponding spaces of equivalence classes are denoted
$\Lr(\ZZ,R)$ and $\Mr(\ZZ,R).$
\end{definitions}

Of course $\Mr\subset \Lr.$ Note that if  $\rho$ doesn't depend on
$z$ and $\rho(s_o)=\infty$ for some $s_o>0,$ $\Mr$ reduces to the
null space and if in addition $R$ is bounded, $\Lr$ is $L_{\infty
}.$ On the other hand, if $\rho $ is a finite function which
doesn't depend on $z$ and $R$ is bounded, $\Mr$ contains all the
bounded functions.

Duality in Orlicz spaces is intimately linked with the convex
conjugacy. The convex conjugate ${\rho^*}$ of $\rho $ is defined
by $\rho^*(z,t)=\sup_{s\in\mathbb{R}}\{st-\rho(z,s)\}.$ It is also
a Young function so that one may consider the Orlicz space $\Lrs.$

\begin{theorem}[Representation of $\Mr'$]\label{res-B1}
Suppose that $\rho$ is a finite Young function. Then, the dual
space of $\Mr$ is isomorphic to $\Lrs.$
\end{theorem}
\proof For a proof of this result, see (\cite{Koz77}, Thm 4.8).
\endproof
 A continuous linear form $\ell\in\Lr'$ is said to be \emph{singular}
 if for all $u\in\Lr,$ there exists a decreasing sequence of
 measurable sets $(A_n)$ such that $R(\cap_n A_n)=0$ and for all
 $n\geq 1,$ $\langle \ell,u\1_{\ZZ\setminus A_n}\rangle=0.$ Let
 us denote $\Lrsi$ the subspace of $\Lr'$ of all singular forms.

\begin{theorem}
 [Representation of $\Lr'$]\label{res-B2}
 Let  $\rho $ be any Young function. The dual space of $\Lr$ is isomorphic to the direct sum
$
 \Lr'=(\Lrs\cdot R)\oplus \Lrsi.
 $
  This implies that any $\ell\in \Lr'$ is uniquely
 decomposed as
 \begin{equation}\label{decomp}
 \ell=\la+\lsing
\end{equation}
with $\la\in\Lrs\cdot R$ and $\lsing\in\Lrsi.$
 \end{theorem}

\proof When $\Lr=L_{\infty}$  this result is the usual
representation of $L_{\infty}'.$
\\
When $\rho$ is a finite function, this result is (\cite{Koz79},
Theorem 2.2).
\\
The general result is proved in \cite{RaoRen}, with $\rho$ not
depending on $z$ but the extension to a $z$-dependent $\rho$ is
obvious.
 \endproof
In the decomposition (\ref{decomp}), $\la$ is called the
\emph{absolutely continuous} part of $\ell$ while $\lsing$ is its
\emph{singular part}.

\begin{proposition}\label{res-B4}
 Let us assume that $\rho$ is finite. Then, $\ell\in\Lr'$ is
 singular if and only if $\langle \ell,u\rangle =0,$ for all $u$ in
 $\Mr.$
 \end{proposition}
 \proof
 This result is (\cite{Koz79}, Proposition 2.1).
 \endproof

The function $\rho $ is said to satisfy the $\Delta_2$-condition
if
\begin{equation}
\label{delta2}
 \textrm{ there exist } C>0,
 s_o\geq 0 \textrm{ such that } \forall s\geq s_o, \rho (2s)\leq
 C\rho (s)
\end{equation}
 If $s_o=0,$ the $\Delta_2$-condition
is said to be global. When $R$ is bounded, in order that
$\Mr=\Lr,$ it is enough that $\rho $ satisfies the
$\Delta_2$-condition. When $R$ is unbounded, this equality still
holds if the $\Delta_2$-condition is global. Consequently, if
$\rho $ satisfies the $\Delta_2$-condition we have $\Lr'=\Lrs\cdot
 R$ so that $\Lrsi$ reduces to the null vector space.

\subsection{The minimization problem ($\pc$)}
Before introducing an extended minimization problem, let us state
properly the basic problem (\ref{eq-01}).

\subsubsection*{Relevant Orlicz spaces}
Since  $\gamma^*_z$ is closed convex for each $z,$ it is the
convex conjugate of some closed convex function $\gamma_z.$
Defining
\begin{equation*}%\label{eq-51}
    \lambda(z,s)=\gamma(z,s)-m(z)s,\quad z\in \ZZ, s\in\mathbb{R}
\end{equation*}
where $m$ satisfies (\ref{eq-04}), one sees that for $R$-a.e.\!
$z,$ $\lambda_z$ is a nonnegative convex function and it vanishes
at 0. Hence,
$$
\lmax(z,s)=\max[\lambda(z,s),\lambda(z,-s)]\in [0,\infty], \quad
z\in \ZZ, s\in \mathbb{R}
$$
is a Young function. We shall use  Orlicz spaces associated with
$\lmax$ and $\lmaxs.$
\\
We denote the space of $R$-absolutely continuous signed measures
having a density in the Orlicz space $\Lls$ by $\LlsR.$ The
effective domain of $I$ is included in $mR+\LlsR.$

\subsubsection*{Constraint}
In order to define the constraint, take $\XXo$ a vector space and
a function $\theta:\ZZ\rightarrow\XXo.$ One wants to give a
meaning to the formal constraint $ \IZ \theta\,dQ=x $ with $Q\in
\LlsR$ and $x\in\XXo.$ Suppose that $\XXo$ is the algebraic dual
space of some vector space $\YYo$ and define for all $y\in\YYo,$
\begin{equation}\label{eq-17}
  \To^\ast y(z):=\langle y,\theta(z)\rangle_{\YYo,\XXo},\
  z\in \ZZ.
\end{equation}
 Assuming that
    \begin{equation}\label{A-bad}
    \To^\ast \YYo\subset\LLl,
    \end{equation}
H\"older's inequality in Orlicz spaces allows to define the
constraint operator $ \To \ell:=\IZ\theta\,d\ell$ for each
$\ell\in\LlsR$ by
\begin{equation}\label{eq-06}
    \left\langle y, \IZ\theta \,d\ell\right\rangle_{\YYo,\XXo} = \IZ \langle
y,\theta(z)\rangle_{\YYo,\XXo}\,\ell(dz),\quad\forall y\in\YYo.
\end{equation}

\subsubsection*{Minimization problem}
Consider the minimization problem
\begin{equation}\label{pco}
    \textsl{minimize } I(Q)
    \textsl{ subject to } \IZ\theta\,d(Q-mR)\in \Co,\quad Q\in mR+\LlsR \tag{$\pco$}
\end{equation}
where $\Co$ is a convex subset of $\XXo.$ One sees with
$\gamma^*_z(t)=\lambda^*_z(t-m(z))$ that
$I_{\gamma^*}(Q)=\Ils(Q-mR).$ Therefore, the problem (\ref{pco})
is equivalent to
\begin{equation}\label{PEnt}
    \textsl{minimize } \Ils(\ell)
    \textsl{ subject to } \IZ\theta\,d\ell\in \Co,\quad \ell\in
    \LlsR
\end{equation}
with $\ell=Q-mR.$ If  the function $m$ satisfies
    $
    m\in\Lls,
    $
one sees with (\ref{A-bad}) and H\"older's inequality in Orlicz
spaces that the vector $x_o=\IZ \theta m\,dR\in\XXo$ is
well-defined in the weak sense. Therefore, (\ref{pco})  is
\begin{equation}\label{pc}
      \textsl{minimize } I(Q)
    \textsl{ subject to } \IZ\theta\,dQ\in C,\quad Q\in \LlsR \tag{$\pc$}
\end{equation}
with $C=x_o+\Co$

\subsection{The extended minimization problem ($\pbc$)}
If the Young function $\lmax$ doesn't satisfy the
$\Delta_2$-condition (\ref{delta2}), for instance if it has an
exponential growth at infinity as in (\ref{eq-02}) or even worse
as in (\ref{eq-03}), the \emph{small} Orlicz space $\EEl$ may be a
proper subset of $\LLl.$ Consequently, for some functions
$\theta,$ the integrability property
\begin{equation}\label{A-good}
    \To^\ast\YYo\subset\EEl
\end{equation}
or equivalently
\begin{equation}\label{A-forall}
 \forall y\in\YYo, \IZ \lambda(\langle y,
\theta\rangle)\, dR<\infty  \tag{A$_\theta^\forall$}
  \end{equation}
may not be satisfied while the weaker property (\ref{A-bad}):
$\To^\ast\YYo\subset\LLl,$ or equivalently
\begin{equation}\label{A-exists}
    \forall y\in\YYo, \exists\alpha>0, \IZ \lambda(\alpha\langle y,
\theta\rangle)\,
  dR<\infty \tag{A$^\exists_\theta$}
\end{equation}
holds. In this situation, analytical complications occur (see
Section \ref{sec-pbc}). This is the reason why constraints
satisfying (\ref{A-forall}) are called \emph{good constraints},
while constraints satisfying (\ref{A-exists}) but not
(\ref{A-forall}) are called \emph{bad constraints}.

If the constraint is bad, it may happen that ($\pc$) is not
attained in $\LlsR.$  This is the reason why it is worth
introducing its  extension ($\pbc$) which may admit minimizers and
is defined by
\begin{equation}\label{pbc}
     \textsl{minimize } \Ib(\ell) \textsl{ subject to }
  \langle\theta,\ell\rangle\in C, \quad\ell\in\Ll'
  \tag{$\pbc$}
\end{equation}
where $\Ll'$ is the topological dual space of $\Ll,$ $\Ib$ and
$\langle\theta,\ell\rangle$ are defined below.
\\
The dual space $\Ll'$ admits the representation $\Ll'\simeq
\LlsR\oplus\Ll^s.$ This means that any $\ell\in\Ll'$ is uniquely
decomposed as $\ell=\ell^a+\ell^s$ where $\ell^a\in\LlsR$ and
$\ell^s\in\Ll^s$ are respectively the \emph{absolutely continuous
part} and the \emph{singular part} of $\ell,$ see Theorem
\ref{res-B2}. The extension $\Ib$ has the following form
\begin{equation}\label{III}
    \Ib(\ell)=I(\ell^a)+\iota_{\dom I_{\gamma}}^*(\ell^s),\quad \ell\in\Ll'
\end{equation}
It will be shown that $\Ib$ is the greatest convex
$\sigma(\Ll',\Ll)$-\lsc\ extension of $I$ to $\Ll'\supset\Lls.$ In
a similar way to (\ref{eq-06}), the assumption
($A_\theta^\exists$) allows to define
$\To\ell=\langle\theta,\ell\rangle$ for all $\ell\in \Ll'$  by
\[
\Big\langle y,\langle\theta,\ell\rangle\Big\rangle_{\YYo,\XXo}
 =\Big\langle\langle y,\theta\rangle ,\ell\Big\rangle_{\Ll,\Ll'},
 \quad \forall y\in\YYo.
\]
Important examples of entropies with $\lmax$ not satisfying the
$\Delta_2$-condition are the usual (Boltzmann) entropy and its
variants, see Section \ref{sec-expl-ent} and (\ref{eq-02}) in
particular.
\\
When $\lmax$ satisfies the $\Delta_2$-condition (\ref{delta2}),
($\pbc$) is ($\pc$).

\subsection{Assumptions}\label{sec-assumptions}
Let us collect the assumptions on $R, \gamma^*$ and $\theta.$

\par\medskip\noindent\textbf{Assumptions (A).}\
\begin{itemize}
  \item[(A$_R$)]
    It is assumed that the reference measure $R$ is a $\sigma$-finite
    nonnegative measure on a space $\ZZ$ endowed with some
    $R$-complete $\sigma$-field.
  \item[(A$_{\gamma^*}$)] \emph{Assumptions on $\gamma^*.$}
     \begin{enumerate}
    \item
    $\gamma^*(\cdot,t)$ is $z$-measurable for all $t$ and
    for $R$-almost every $z\in \ZZ,$ $\gamma^*(z,\cdot)$
    is a \lsc\ strictly convex $[0,+\infty]$-valued function on
    $\mathbb{R}$ which attains its (unique) minimum at $m(z)$ with
    $\gamma^*(z,m(z))=0.$
    \item
    $\IZ \lambda^*(\alpha m)\,dR+\IZ \lambda^*(-\alpha
    m)\,dR<\infty,$ for some $\alpha>0.$
     \end{enumerate}
  \item[(A$_\theta$)] \emph{Assumptions on $\theta.$}
    \begin{enumerate}
    \item for any $y\in\YYo,$ the function $z\in\ZZ\mapsto\langle
    y,\theta(z)\rangle\in\R$ is measurable;
     \item for any $y\in\YYo,$ $\langle y,\theta(\cdot)\rangle=0, R\hbox{-a.e.}$ implies that $y=0;$
     \item[($\exists$)] $\forall y\in\YYo, \exists\alpha>0,\quad \IZ \lambda(\alpha\langle y, \theta\rangle)\,
        dR<\infty.$
    \end{enumerate}
\end{itemize}

\begin{remarks}\label{rem-02}
Some technical remarks about the assumptions.
    \begin{enumerate}[(a)]
    \item
    Since $\gamma_z^*$ is a convex function on $\mathbb{R},$  it is continuous on the interior of its
    domain. Under our assumptions, $\gamma^*$ is (jointly) measurable,
    and so are $\gamma$ and $m.$ Hence,  $\lambda$ is also
    measurable.
    \item
    As $\gamma_z^*$ is strictly convex, $\gamma_z$ is differentiable.
    \item
    Assumption (A$_{\gamma^*}^2$) is $m\in\Lls.$ It allows to
    consider Problem (\ref{pc}) rather than (\ref{pco}). If this assumption
    is not satisfied, our results still hold for (\ref{pco}), but their statement is a little heavier, see Remark
    \ref{rem-01}-d below.
    \item
    Since $\XXo$ and $\YYo$ are in separating
    duality, ($A_\theta^2$) states that the
    vector space spanned by the range of $\theta$ ``is essentially"
    $\XXo.$ This is not an effective restriction.
    \end{enumerate}
\end{remarks}

\section{Preliminary results}\label{sec-prelim}

\newcommand{\Ci}{C\cap\XX}

The aim of this section is to recall for the convenience of the
reader some results of  \cite{Leo01a, Leo03, p-Leo07a}.

\subsection{Convex minimization problems under weak constraint qualifications}\label{sec-saddle}
The main results of \cite{p-Leo07a} are presented.

\subsubsection*{Basic diagram}%\label{sec:Tadj}
Let $\UUo$ be a vector space, $\LLo=\UUo^\ast$ its algebraic dual
space, $\Fo$ a $(-\infty,+\infty]$-valued convex function on
$\UUo$ and $\Fos$ its convex conjugate for the duality $\langle
\UUo,\LLo\rangle:$
\begin{equation*}
    \Fos(\ell):= \sup_{u\in\UUo}\{\ul-\Fo(u)\},\quad \ell\in\LLo\\
\end{equation*}
 Let $\YYo$ be another vector space,
$\XXo=\YYo^\ast$ its algebraic dual space and
$\To:\LLo\rightarrow\XXo$ a linear operator. We consider the
convex minimization problem
\begin{equation}
  \textsl{minimize } \Fos(\ell) \textsl{ subject to } \To\ell\in  C,\
  \ell\in\LLo
\tag{$\Po$}
\end{equation}
where $ C$ is a convex subset of $\XXo.$
\\
This will be used later with $\Fi=\Il$ on the Orlicz space
$\UUo=\mathcal{E}_\lmax(\ZZ,R)$ or
$\UUo=\mathcal{L}_\lmax(\ZZ,R).$
\\
It is useful to define the constraint operator $\To$ by means of
its adjoint $\To^\ast :\YYo\rightarrow\LLo^\ast$ for each
$\ell\in\LLo,$ by
    $
  \langle\To^\ast y,\ell\rangle_{\LLo^\ast,\LLo}= \langle
  y,\To\ell\rangle_{\YYo,\XXo},$ $
\forall y\in\YYo.
    $

\subsubsection*{Hypotheses}

Let us give the list of the main hypotheses.
\begin{itemize}
\item[\HF]
1-\quad $\Fo: \UUo\rightarrow [0,+\infty]$ is $\sUL$-\lsc, convex and $\Fo(0)=0$\\
2-\quad $\forall u\in\UUo, \exists \alpha >0, \Fo(\alpha u)<\infty$\\
3-\quad $\forall u\in\UUo, u\not=0,\exists t\in\R, \Fo(tu)>0$
\item[\HT]
1-\quad $\To^\ast (\YYo)\subset\UUo$\\
2-\quad $\mathrm{ker\ }\To^\ast =\{0\}$ \item [\HC] \qquad $
C\cap\XX$ is a convex $\sigma(\XX,\YY)$-closed subset of $\XX$
\end{itemize}
The definitions of the vector spaces $\XX$ and $\YY$ which appear
in the last assumption are stated below. For the moment, let us
only say that if $  C$ is convex and $\sigma(\XXo,\YYo)$-closed,
then \HC\ holds.

\subsubsection*{Several primal and dual problems.}\quad %\label{sec:problems}
 These variants are expressed below in terms of new spaces and
functions. Let us first introduce them.

\begin{enumerate}[-]
\item \textsf{The norms $|\cdot|_\Fi$ and $|\cdot|_\La.$}\quad Let
$\Fi_{\pm}(u)=\max(\Fo(u),\Fo(-u)).$ By \HFi\ and \HFii,
$\{u\in\UUo; \Fi_{\pm}(u)\leq 1\}$  is a convex absorbing balanced
set. Hence its gauge functional which is defined for all $
u\in\UUo$ by $|u|_\Fi := \inf\{\alpha>0; \Fi_{\pm}(u/\alpha))\leq
1\}$ is a seminorm. Thanks to hypothesis \HFiii, it is a norm.
\\
Taking \HTi\ into account, one can define
\begin{equation}\label{La}
  \Lao(y):= \Fo(\To^\ast y), y\in\YYo.
\end{equation}
Let $\La_\pm(y)= \max(\Lao(y),\Lao(-y)).$ The gauge functional on
$\YYo$ of the set $\{y\in\YYo; \La_\pm(y)\leq 1\}$ is $|y|_\La :=
\inf\{\alpha>0;\La_\pm(y/\alpha)\leq 1\}, y\in\YYo.$ Thanks to
\HF\ and \HT, it is a norm and
\begin{equation*}%\label{eq-117}
    |y|_\La=|\To^\ast y|_\Fi, \quad y\in\YYo.
\end{equation*}

\item \textsf{The spaces.}\quad Let
\begin{eqnarray*}
& &\UU \mathrm{\ be\ the\ } \NF\textrm{-completion\ of\ } \UUo\mathrm{\ and\ let }\\
  & &\LL:= (\UUo,\NF)'\textrm{\ be\ the\ topological\ dual\ space\ of\ }
  (\UUo,\NF).
\end{eqnarray*}
Of course, we have $(\UU,\NF)'\cong\LL\subset\LLo$ where any
$\ell$
 in $\UU'$ is identified with its restriction to $\UUo.$ Similarly, we introduce
\begin{eqnarray*}
& &\YY \mathrm{\ the\ } \NL\textrm{-completion\ of\ } \YYo
  \mathrm{\ and }\\
 & &\XX:= (\YYo,\NL)'\mathrm{\ the\ topological\ dual\ space\ of\
 }(\YYo,\NL).
\end{eqnarray*}
We have $(\YY,\NL)'\cong\XX\subset\XXo$ where any $x$ in $\YY'$ is
identified with its restriction to $\YYo.$
\\
We also have to consider the \emph{algebraic} dual spaces
$\LL^\ast$ and $\XX^\ast$  of $\LL$ and $\XX.$

\item \textsf{The operators $T$ and $T^\ast.$}\quad Let us denote
$T$ the restriction of $\To$ to $\LL\subset\LLo.$
 One can show that under ($H_{\Fi\& T}$),
$
    \To\LL\subset\XX.
$ Hence $T: \LL\to\XX.$ Let us define its adjoint $T^\ast
:\XX^\ast\rightarrow\LL^\ast$ for all $\omega\in\XX^\ast$ by: $
  \langle  \ell,T^\ast\omega\rangle_{\LL,\LL^\ast}=\langle T\ell,\omega\rangle_{\XX,\XX^\ast},
  \forall \ell\in\LL.
$ We have the inclusions $\YYo\subset\YY\subset\XX^\ast.$ The
adjoint operator $\To^*$ is the restriction of $T^*$ to $\YYo.$

\item \textsf{The functionals.}\quad They are:
$$
\begin{array}{rcll}
    \Fb(\zeta)&:=&\sup_{\ell\in\LL}\{\langle\zeta,\ell\rangle-\Fs(\ell)\},&\quad\zeta\in\LL^*\\
       \La(y)&:=& \Fb(T^*y),&\quad  y\in\YY\\
   \Lb(\omega)&:=&\Fb(T^*\omega),&\quad \omega\in\XX^*\\
     \Los(x)&:=& \sup_{y\in\YYo}\{\yx-\Lao(y)\},&\quad x\in\XXo\\
    \Ls(x)&:=& \sup_{y\in\YY}\{\yx-\La(y)\},&\quad x\in\XX\\
\end{array}
$$
\item \textsf{The optimization problems.}\quad They are:
\begin{align}
&\textsl{minimize } \Fos(\ell) & &\textsl{subject to } \To\ell\in
 C,&& \ell\in\LLo   \tag{$\Po$}\\
&\textsl{minimize } \Fos(\ell) & &\textsl{subject to }
T\ell\in C,&& \ell\in\LL   \tag{$\mathcal{P}$}\\
&\textsl{maximize } \inf_{x\in\Ci}\yx - \La(y), && &
&y\in\YY  & &   \tag{$\mathcal{D}$}\\
&\textsl{maximize } \inf_{x\in\Ci}\xo - \Lb(\omega), && &
&\omega\in\XX^* & & \tag{$\Db$}
\end{align}
\end{enumerate}

\subsubsection*{Statement of the results.} It is assumed that \HF,
\HT\ and \HC\ hold.

\begin{theorem}[Primal attainment and dual equality]\label{xP3}\
\begin{enumerate}[(a)]
    \item The problems $(\Po)$ and $(\mathcal{P})$ are equivalent: they have
        the same solutions and $\inf(\Po)=\inf(\mathcal{P})\in[0,\infty].$

    \item  We have the dual equalities
\begin{equation*}%\label{xed1}
   \inf(\Po)=\inf(\mathcal{P})=\sup(\mathcal{D})=\sup(\Db)=\inf_{x\in  C}\Los(x)=\inf_{x\in \Ci}\Ls(x)\in [0,\infty]
\end{equation*}
    \item If in addition $\{\ell\in\LLo;\To\ell\in
C\}\cap\dom\Fos\not=\emptyset,$ then $(\Po)$ is attained in
    $\LL.$ Moreover, any minimizing sequence for $(\Po)$ has $\sLUi$-cluster
    points and every such cluster point solves $(\Po)$.
\end{enumerate}
\end{theorem}

\begin{theorem}[Dual attainment and representation. Interior convex
constraint]\label{T3a}\ \\ Assume that
$C\cap\icor(\To\dom\Fos)\not=\emptyset.$
\\
Then, the primal problem $(\Po)$ is attained in $\LL$ and the
extended dual problem $(\Db)$ is attained in $\XX^*.$ Any solution
$\lh\in\LL$ of $(\Po)$ is characterized by the existence of some
$\ob\in\XX^*$ such that
 \begin{equation}\label{xeq-96}
    \left\{\begin{array}{cl}
      (a) & T\lh\in  C \\
      (b) & \langle T^*\ob,\lh\rangle\leq \langle T^*\ob,\ell\rangle
      \textrm{ for all }\ell\in\{\ell\in\LL; T\ell\in C\} \cap\dom\Fs\\
      (c) & \lh\in\partial_{\LL}\Fb(T^*\ob) \\
    \end{array}\right.
\end{equation}
Moreover, $\lh\in\LL$ and $\ob\in\XX^*$ satisfy (\ref{xeq-96}) if
and only if $\lh$ solves $(\Po)$ and $\ob$ solves $(\Db)$.
\end{theorem}
The assumption $C\cap\icor(\To\dom\Fos)\not=\emptyset$ is
equivalent to $ C\cap\icordom\Los\not=\emptyset$ and the
representation formula (\ref{xeq-96}-c)
 is equivalent to Young's identity
\begin{equation}\label{xeq-11}
    \Fs(\lh)+\Fb(T^*\ob)=\langle \ob,T\lh\rangle=\Ls(\xh)+\Lb(\ob).
\end{equation}
Formula (\ref{xeq-96}-c) can be made a little more precise by
means of the following regularity result.

\begin{theorem}\label{res-01}
Any solution $\ob$ of $(\Db)$ shares the following properties
\begin{itemize}
 \item[(a)] $\ob$ is in the $\sigma(\XX^*,\XX)$-closure of $\dom\La;$
 \item[(b)] $T^\ast \ob$ is in the $\sigma(\LL^*,\LL)$-closure of
 $T^\ast(\dom\La).$
\end{itemize}
If in addition the level sets of $\Fi$ are $\NF$-bounded, then
\begin{itemize}
 \item[(a')] $\ob$ is in $\YY''.$ More precisely, it is in the $\sigma(\YY'',\XX)$-closure of $\dom\La;$
 \item[(b')] $T^\ast \ob$ is in $\UU''.$ More precisely, it is in the $\sigma(\UU'',\LL)$-closure of
 $T^\ast(\dom\La)$
\end{itemize}
where $\YY''$ and $\UU''$ are the \emph{topological} bidual spaces
of $\YY$ and $\UU.$ This occurs if $\Fi,$ or equivalently $\Fs,$
is an even function.
\end{theorem}

\subsection{Convex conjugates in a Riesz space}\label{sec-Riesz}

The following results are taken from \cite{Leo01a, Leo03}. For the
basic definitions and properties of Riesz spaces, see
\cite[Chapter 2]{Bou-Int1-4}.
\\
Let $U$ be a Riesz vector space for the order relation $\leq .$
Since $U$ is a Riesz space, any $u\in U$ admits a nonnegative
part: $ u_+:=u\vee 0,$ and a nonpositive part: $ u_-:= (-u)\vee
0.$ Of course, $u= u_+- u_-$ and as usual, we state: $\vert u\vert
= u_++ u_-.$

\begin{remark}\label{rem-04}
Recall that there is a natural order on the algebraic dual
 space $E^*$ of a Riesz vector space $E$ which is defined by:
 $e^*\leq f^*\Leftrightarrow\langle e^*,e\rangle \leq \langle
 f^*,e\rangle $ for any $e\in E$ with $e\geq 0.$ A linear form
 $e^*\in E^*$ is said to be \emph{relatively bounded} if
 for any $f\in E,$ $f\geq 0,$ we have $\sup_{e:\vert e\vert
 \leq f} \vert \langle e^*,e\rangle \vert <+\infty .$ Although
 $E^*$ may not be a Riesz space in general, the vector space
 $E^b$ of all the relatively bounded linear forms on $E$ is always
 a Riesz space. In particular, the elements of $E^b$ admit a
 decomposition in positive and negative parts
 $e^*=e^*_+-e^*_-.$
\end{remark}

Let $\Phi $ be a $[0,\infty ]$-valued function on $U$ which
satisfies the following conditions:
\begin{align}
    &\forall u\in U,\ \Phi (u)=\Phi ( u_+- u_-)=\Phi( u_+)+\Phi (- u_-)\label{eq-30a}\\
    &\forall u,v\in U,\ \left\{ \begin{array}{ccc}
        0\leq u\leq v & \Longrightarrow & \Phi (u)\leq \Phi (v)\\
        u\leq v\leq 0 & \Longrightarrow & \Phi (u)\geq \Phi (v)\\
        \end{array}\right.\label{eq-30b}
\end{align}
Clearly (\ref{eq-30a}) implies $\Phi (0)=0,$ (\ref{eq-30a}) and
(\ref{eq-30b}) imply that for any $u\in U,$ $\Phi (u)=\Phi (
u_+)+\Phi (- u_-)\geq \Phi (0)+\Phi (0)=0.$ Therefore, $\Fs$ is
$[0,\infty ]$-valued and $\Fs(0)=0.$

For all $u\in U,$ $\Phi_+(u)=\Phi(|u|),$ $\Phi_-(u)=\Phi(-|u|).$
The convex conjugates of $\Phi, \Phi_+$ and $\Phi_-$ with respect
to $\langle U,U^*\rangle$ are denoted $\Fs, \Fs_+$ and $\Fs_-.$
Let $L$ be the vector space spanned by $\dom\Fs.$ The convex
conjugates of $\Fs, \Fs_+$ and $\Fs_-$ with respect to $\langle
L,L^*\rangle$ are denoted $\Fb, \overline{\Phi_+}$ and
$\overline{\Phi_-}.$ The space of relatively bounded linear forms
on $U$ and $L$ are denoted by $U^b$ and $L^b,$ whenever $L$ is a
Riesz space.

One writes $a_\pm\in A_\pm$ for [$a_+\in A_+$ and $a_-\in A_-$].

 \begin{proposition}\label{res-05}
Assume (\ref{eq-30a}) and (\ref{eq-30b}) and suppose that $L$ is a
Riesz space.
\begin{enumerate}[(a)]
    \item For all $\ell\in U^*,$
$$
\Fs( \ell)=\left \{
    \begin{array}{ll}
       \Fs_+(\ell_+)+\Fs_-(\ell_-) &\hbox{ if } \ell\in U^b\\
      +\infty &\hbox{ otherwise} \\
    \end{array}
    \right.
$$
    \item
Denoting $L_+$ and $L_-$ the vector subspaces of $L$ spanned by
$\dom\Fs_+$ and $\dom\Fs_-,$ we have
$$
\Fb(\zeta )=\left \{
    \begin{array}{ll}
      \overline{\Phi_+}({\zeta_+}_{|L_+})+ \overline{\Phi_-}({\zeta_-}_{|L_-} ) &\hbox{ if }\zeta\in L^b \\
      +\infty &\hbox{ otherwise} \\
    \end{array}
\right .
$$
which means that
$\overline{\Phi_\pm}(\zeta_\pm)=\overline{\Phi_\pm}(\zeta_\pm')$
if $\zeta_\pm$ and $\zeta_\pm'$ match on $L_\pm.$
    \item
     Let $\ell\in L,$ $\zeta\in L^*$ be such that
$\ell\in\partial_L\Fb(\zeta).$ Then, $\ell_\pm\in
\partial_{L_\pm} \overline{\Phi_\pm}({\zeta_\pm}_{|L_\pm})\subset L_\pm.$
\end{enumerate}
\end{proposition}

\proof (a) and (b) are proved at \cite[Proposition 4.4]{Leo01a}
under the additional assumption that for all $u\in U$ there exists
$\lambda >0$ such that $\Phi(\lambda u)<+\infty.$ But it can be
removed. Indeed, if for instance $\Phi_-$ is null, $\Fs_-$ is the
convex indicator of $\{0\}$ whose domain is in $U^b.$ The
statement about $\Fb$ is an iteration of this argument.
\\
The last statement of (b) about ${\zeta_\pm}_{|L_\pm}$ directly
follows from $\dom\Fs_\pm\subset L_\pm.$
\\
For (c), see the proof of \cite[Proposition 4.5]{Leo03}.
\endproof

\section{Solving ($\pbc$)}\label{sec-pbc}

The general assumptions (A) are imposed and we study (\ref{pbc}).

\subsection{Several function spaces and cones}
To  state the extended dual problem (\ref{dbc}) below, notation is
needed. If $\lambda$ is not an even function, one has to consider
\begin{equation}\label{eq-14}
    \left\{
    \begin{array}{l}
       \lambda_+(z,s)=\lambda(z,|s|) \\
       \lambda_-(z,s)=\lambda(z,-|s|) \\
    \end{array}
    \right.
\end{equation}
which are Young functions and the corresponding Orlicz spaces.

\begin{definitions}\label{def-1}
For any relatively
 bounded linear form $\zeta$ on $\Ll'$ i.e.\! $\zeta\in\Ll'^b,$ one writes:
\begin{itemize}
    \item  $\zeta\in K_{\lambda}''$ to specify that
    ${\zeta_\pm}_{|L_{\lambda_\pm}'\cap\Ll'}\in L_{\lambda_\pm}''$
    \item  $\zeta\in K_{\lambda^*}'$ to specify that
    ${\zeta_\pm}_{|L_{\lambda_\pm^*}R\cap\Ll'}\in L_{\lambda_\pm^*}'$
    \item $\zeta\in K_{\lambda}$ to specify that
    ${\zeta_\pm}_{|L_{\lambda_\pm^*}R\cap\Ll'}\in L_{\lambda_\pm}$
    \item $\zeta\in K_{\lambda^*}^s$ to specify that
    ${\zeta_\pm}_{|L_{\lambda_\pm^*}R\cap\Ll'}\in L_{\lambda_\pm^*}^s$
    \item $\zeta\in K_{\lambda}^{s\prime}$ to specify that
    ${\zeta_\pm}_{|L_{\lambda_\pm}^s\cap\Ll'}\in L_{\lambda_\pm}^{s\prime}$
\end{itemize}
where $\lambda_\pm$ are defined at (\ref{eq-14}) and
${\zeta_\pm}_{|L_\pm\cap\Ll'}\in L_\pm'$ means that the
restriction of $\zeta_\pm$ to $L_\pm\cap \Ll'$ is continuous with
respect to relative topology generated by the strong topology of
$L_\pm$ on $L_\pm\cap \Ll'.$
\begin{enumerate}
    \item The sets $K_{\lambda}'', K_{\lambda^*}', K_{\lambda},
K_{\lambda^*}^s$ and $K_{\lambda}^{s\prime}$ are defined to be the
corresponding subsets of $\Ll'^b.$ They  are not vector spaces in
general but convex cones with vertex 0.
    \item The $\sigma(K_{\lambda}'',K_{\lambda}')$-closure $\overline{A}$ of
a set $A$ is defined as follows: $\zeta\in\Ll^{\prime b}$ is in
$\overline{A}$ if ${\zeta_\pm}_{|L_{\lambda_\pm}'\cap\Ll'}$ is in
the
$\sigma(L_{\lambda_\pm}''\cap\Ll'^b,L_{\lambda_\pm}'\cap\Ll')$-closure
of $A_\pm=\{\zeta_\pm;\zeta\in A\}.$ Clearly,
$\overline{A}_\pm=\{\zeta_\pm;\zeta\in\overline{A}\}.$
    \\ One defines similarly the  $\sigma(K_{\lambda^*}',K_{\lambda^*}),$
$\sigma(K_{\lambda},K_{\lambda}'),$
$\sigma(K_{\lambda^*}^s,K_{\lambda^*})$ and
$\sigma(K_{\lambda}^{s\prime},K_{\lambda}^{s})$-closures.
    \item Let $A$ be a subset of $\Ll.$ Its strong closure $s\textrm{-}\cl A$ in
    $K_\lambda$ is the set of all measurable functions $u$ such
    that $u_\pm$ is in the  $\|\cdot\|_{\lambda_\pm}$-closure of $A_\pm=\{v_\pm;v\in A\}.$
\end{enumerate}
\end{definitions}

Let $\rho$ be a Young function. By Theorem \ref{res-B2}, we have
$L_\rho''=[L_{\rho}. R\oplus L_{\rho}^s]\oplus
L_{\rho}^{s\prime}.$ For any $\zeta\in L_{\rho}''=(
L_{\rho^*}R\oplus L_{\rho}^s)',$ let us denote the restrictions
$\zeta_1=\zeta_{| L_{\rho^*}R}$ and $\zeta_2=\zeta_{|
L_{\rho}^s}.$ Since, $( L_{\rho^*}R)'\simeq L_{\rho}\oplus
L_{\rho^*}^s,$ one sees that any $\zeta\in L_{\rho}''$ is uniquely
decomposed into
\begin{equation}\label{eq-15}
    \zeta=\zeta_1^a+\zeta_1^s+\zeta_2
\end{equation}
with $\zeta_1=\zeta_1^a+\zeta_1^s\in L_{\rho^*}',$ $\zeta_1^a\in
L_{\rho},$ $\zeta_1^s\in L_{\rho^*}^s$ and $\zeta_2\in
L_{\rho}^{s\prime}.$ With our definitions,
$K_\lambda''=[K_{\lambda}\oplus K_{\lambda^*}^s]\oplus
K_{\lambda}^{s\prime}$ and the decomposition (\ref{eq-15}) holds
for any $\zeta\in K_\lambda''$ with
\begin{equation*}
    \left\{
\begin{array}{l}
  \zeta_1=\zeta_1^a+\zeta_1^s\in K_{\lambda}\oplus K_{\lambda^*}^s=K_{\lambda^*}',\\
   \zeta_2\in K_{\lambda}^{s\prime}. \\
\end{array}%
\right.
\end{equation*}

\subsection{The ingredients of the saddle-point method}
One applies the abstract results of Section \ref{sec-saddle} with
\begin{equation}\label{eq-09}
     \Fo(u)=\Il(u):= \IZ\lambda(u)\,dR,\quad u\in\UUo:=\LLl
\end{equation}
This gives $\UU=\Ll$ with the Orlicz norm $|u|_\Fi=\|u\|_\lmax$
and $\LL=\Ll'=\LlsR\oplus\Ll^s,$ by Theorem \ref{res-B2}. The
space $\YY$ is the completion of $\YYo$ endowed with the norm
$|y|_\Lambda=\|\yt\|_\lmax.$ One denotes $\YY=\YY_L.$ It is
isomorphic to the closure of the subspace $\{\yt;y\in\YYo\}$ in
$\Ll,$ see assumption (\ref{A-exists}). With some abuse of
notation, one still denotes $T^*y=\yt$ for $y\in\YY_L.$ Remark
that this can be interpreted as a dual bracket between $\XXo^*$
and $\XXo$ since $T^*y=\langle\tilde{y},\theta\rangle$ $R$-a.e.\!
for some $\tilde{y}\in\XXo^*.$ The topological dual space
$\XX_L=\YY_L'$ is identified with $\Ll'/\mathrm{ker\,}T$
 and its norm is given by $|x|_\Lambda^*=\inf\{\|\ell\|_\lmax^*;\ell\in\Ll':
 T(\ell)=x\}.$ This last identity is a dual equality as in Theorem
 \ref{xP3}-b with $\Phi=\iota_B$ where $B$ is the unit ball of
 $\Ll$ and $C=\{x\}.$
 \\
 The assumption ($H_C$) that $C$ is $\sigma(\XX_L,\YY_L)$-closed
 convex is equivalent to
 \begin{equation}\label{eq-10}
   \To^{-1}C\cap\Ll'=\bigcap_{y\in Y}\left\{\ell\in\Ll';\langle\yt,\ell\rangle \ge a_y\right\}
\end{equation}
for some subset $Y\subset\YY_L$ and some functions $y\in Y\mapsto
a_y\in\R.$ For comparison, note that if $C$ is only supposed to be
convex, $\bigcap_{(y,a)\in
A}\left\{\ell\in\Ll';\langle\yt,\ell\rangle > a\right\}$ with
$A\subset \YY\times\R$ is the general shape of $T^{-1}C.$

\subsection{The main result}
Let us define
\begin{equation*}
  \Gamma^*(x) =
  \sup_{y\in\YYo}\left\{\yx-I_\gamma(\yt)\right\},\quad x\in\XXo
\end{equation*}
which is the convex conjugate of
    $
\Gamma(y) =I_\gamma(\yt),$ $y\in\YYo.
    $
The dual problem ($\mathcal{D}$) associated with (\ref{pc}) and
(\ref{pbc}) is
\begin{equation}\label{dc}
    \textsl{maximize } \inf_{x\in \Ci}\yx-I_\gamma(\yt),
   \quad y\in\YY\tag{$\dc$}
\end{equation}
The extended dual problem is
\begin{equation}\label{dbc}
   \textsl{maximize } \inf_{x\in C\cap\XX_L}\langle\omega,x\rangle-
   \Il\big([T^*\omega]_1^a\big)+\iota^*_{\dom\Ils}\big([T^*\omega]_1^s\big)+\iota_D\big([T^*\omega]_2\big),
   \quad \omega\in\YYb \tag{$\dbc$}
\end{equation}
where
\begin{itemize}
    \item $T^*:\XX_L^*\to \Ll'^*$ is the extension of $\To^*$ which is
defined at Section \ref{sec-saddle},
    \item $D$ is the $\sigma(K_{\lambda}^{s\prime},K_{\lambda}^s)$-closure of $\dom
    I_{\lambda}$ and
    \item $\YYb$ is the cone of all $\omega\in\XX_L^*$ such that $T^*\omega\in K_\lambda''.$
\end{itemize}
Clearly,
    $
  \iota^*_{\dom\Ils}(\zeta_1^s)=\iota^*_{\dom I_{\lambda_+^*}}(\zeta_{1+}^s)+\iota^*_{\dom
    I_{\lambda_-^*}}(\zeta_{1-}^s)$ and
    $
  \iota_D(\zeta_2)=\iota_{D_+}(\zeta_{2+})+\iota_{D_-}(\zeta_{2-})
    $
    where
    $D_\pm$ is the $\sigma(L_{\lambda_\pm}^{s\prime}\cap\Ll'^b,L_{\lambda_\pm}^s\cap\Ll')$-closure of $\dom
    I_{\lambda_\pm}.$
\\
As $R$ is assumed to be $\sigma$-finite, there exists a measurable
partition $\seq {\ZZ}k$ of $\ZZ:$ $\bigsqcup_k\ZZ_k=\ZZ,$ such
that $R(\ZZ_k)<\infty$ for each $k\ge1.$
\begin{theorem}\label{res-03}
Suppose that
\begin{enumerate}
    \item the assumptions \emph{(A)} are satisfied;
    \item for each $k\ge1,$ $\Ll(\ZZ_k,R_{|\ZZ_k})$ is dense in $L_{\lambda_+}(\ZZ_k,R_{|\ZZ_k})$  and
    $L_{\lambda_-}(\ZZ_k,R_{|\ZZ_k})$ with respect to the topologies associated with  $\|\cdot\|_{\lambda_+}$ and
    $\|\cdot\|_{\lambda_-};$
    \item $C$ satisfies (\ref{eq-10}) with $\yt\in\Ll$ for all $y\in Y.$
\end{enumerate}
Then:
\begin{enumerate}[]
\item[(a)] The dual equality for \emph{(\ref{pbc})} is
\begin{equation*}
    \inf(\pbc)=\inf_{x\in C}\Gamma^*(x)=\sup(\dc)=\sup(\dbc)\in [0,\infty].
    \end{equation*}
 \item[(b)] If $C\cap \dom\Gamma^*\not=\emptyset$ or equivalently $C\cap \To\dom \Ib\not=\emptyset,$
 then \emph{(\ref{pbc})}
 admits  solutions in $\Ll',$ any minimizing sequence admits $\sigma(\Ll',\Ll)$-cluster points
 and every such point is a solution to \emph{(\ref{pbc})}.
\end{enumerate}
 Suppose that in addition we have
 \begin{equation}\label{eq-24}
    C\cap \icordom\Gamma^*\not=\emptyset
\end{equation}
 or equivalently $C\cap \icor(\To\dom \Ib)\not=\emptyset.$ Then:
\begin{enumerate}[]
 \item[(c)]
 Let us denote $\xh\eqdef T^*\lh.$ There exists  $\ob\in\YYb$ such that
 \begin{equation}\label{eq-22}
    \left\{\begin{array}{cl}
      (a) & \xh\in C\cap\dom\Gamma^* \\
      (b) & \langle \ob,\xh\rangle_{\XX_L^*,\XX_L} \leq \langle \ob,x\rangle_{\XX_L^*,\XX_L}, \forall x\in C\cap\dom\Gamma^* \\
      (c) & \lh\in\gamma'_z([T^*\ob]_1^a)\,R+D^\bot([T^*\ob]_2) \\
    \end{array}\right.
\end{equation}
where
$$
D^\bot(\eta)=\{k\in \Ll^s; \forall h\in \Ll, \eta+h\in D
\Rightarrow \langle h,k\rangle\le0\}
$$
is the outer normal cone of $D$ at $\eta.$
\\
$T^*\ob$ is in the $\sigma(K_\lambda'',K_\lambda')$-closure of
$T^*(\dom\Lambda)$ and there exists some $\ot\in\XXo^*$ such that
$$
[T^*\ob]_1^a=\langle\ot,\theta(\cdot)\rangle_{\XXo^*,\XXo}
$$
is a measurable function in the strong closure of
$T^*(\dom\Lambda)$ in $K_\lambda.$
\\
Furthermore, $\lh\in\Ll'$ and $\ob\in\YYb$ satisfy (\ref{eq-22})
if and only if $\lh$ solves \emph{(\ref{pbc})} and $\ob$ solves
\emph{(\ref{dbc})}.
    \item[(d)]
     Of course, (\ref{eq-22}-c) implies
    $
\xh=\IZ \theta\gamma'(\ott)\,dR+\langle\theta,\lh^s\rangle.
    $
    Moreover,
    \begin{enumerate}[1.]
        \item $\xh$ minimizes $\Gamma^*$ on $C,$
        \item
        $\Ib(\lh)=\Gamma^*(\xh)=\IZ\gamma^*\circ\gamma'(\ott)\,dR+\sup\{\langle u,\lh^s\rangle;u\in\dom I_\gamma\}<\infty$
        and
        \item $\Ib(\lh)+\IZ\gamma(\ott)\,dR=\IZ\ott\,d\lh^a+\langle [T^*\ob]_2,\lh^s \rangle_{{K_\lambda^s}',K_\lambda^s}.$
    \end{enumerate}
\end{enumerate}
\end{theorem}

\begin{proposition}\label{res-06}
For the assumption (2) of Theorem \ref{res-03} to be satisfied, it
is enough that one of these conditions holds
    \begin{enumerate}[(i)]
        \item
        $\lambda$ is even or more generally
        $0<\Liminf t \frac{\lambda_+}{\lambda_-}(t)\le\Limsup t \frac{\lambda_+}{\lambda_-}(t)<+\infty;$
        \item
        $\Lim t \frac{\lambda_+}{\lambda_-}(t)=+\infty$ and $\lambda_-$
        satisfies the $\Delta_2$-condition (\ref{delta2}).
    \end{enumerate}
\end{proposition}
\begin{proof}
It is enough to work with a bounded measure $R.$
\\
Condition (i) is equivalent to $L_{\lambda_+}=L_{\lambda_-}=\Ll$
and the result follows immediately.
\\
Condition (ii) says that $\lambda_+=\lmax$ and
$L_{\lambda_-}=E_{\lambda_-}.$ As $\gamma^*$ is assumed to be
strictly convex, zero is in the interior of $\dom\lambda$ and
$\Ll$ contains the space $B$ of all bounded measurable functions.
But $B$  is dense in $E_{\lambda_-}$ and the result follows.
\end{proof}

\begin{remarks}\label{rem-01} General remarks about Theorem
\ref{res-03}.\
\begin{enumerate}[(a)]
      \item
The assumption (3) is equivalent to $C$ is
$\sigma(\XX_L,\YY_L)$-closed convex.
    \item
The dual equality with $C=\{x\}$ gives for all $x\in\XXo$
    \begin{equation*}
  \Gamma^*(x)  =\inf\left\{\Ib(\ell) ;
  \ell\in\Ll',\langle\theta,\ell\rangle=x\right\}.
    \end{equation*}

    \item
Note that $\ob$ does not necessarily belong to $\YYo.$ Therefore,
the Young equality $\langle\ob,\xh\rangle
=\Gamma^*(\xh)+\Gamma(\ob)$ is meaningless. Nevertheless, there
exists a natural extension $\overline{\Gamma}$ of $\Gamma$ such
that $\langle\xh,\ob\rangle =\Gamma^*(\xh)+\overline{\Gamma}(\ob)$
holds, see (\ref{xeq-11}). This gives the statement (d-3).

    \item
Removing the assumption (A$_{\gamma^*}^2$): $m\in\Lls,$ one can
still consider the minimization problem
\begin{equation}\label{pbco}
    \textsl{minimize } \Ib(\ell)
    \textsl{ subject to } \langle\theta,\ell-mR\rangle\in \Co,\quad \ell\in mR+\Ll' \tag{$\pbco$}
\end{equation}
instead of (\ref{pbc}). The transcription of Theorem \ref{res-03}
is as follows.  Denote
\begin{equation*}
  \Lambda^*(x) =
  \sup_{y\in\YYo}\left\{\yx-\IZ\lambda(\yt)\,dR\right\},\quad
  x\in\XXo
\end{equation*}
and replace respectively (\ref{pbc}), $C,$ $\Gamma^*,$ $\xh$ and
$\gamma$ by (\ref{pbco}), $\Co,$ $\Ls,$ $\tilde{x}$ and $\lambda$
where $\tilde{x}=\langle\theta,\lh-mR\rangle$ is well-defined.
\\
The statement (b) must be replaced by the following one:
\textit{If $\Co\cap \dom\Lambda^*\not=\emptyset,$ then
\emph{(\ref{pbco})}
 admits  solutions in $mR+\Ll',$  any minimizing sequence $\seq
    {\ell}n$ is such that $(\ell_n-mR)_{n\ge1}$
    admits cluster points $\lh-mR$ in $\Lls'$ with respect to the topology $\sigma(\Ll',\Ll)$
    and $\lh$ is a solution of \emph{(\ref{pbco})}.}
\end{enumerate}
\end{remarks}

\begin{proof}[Proof of Theorem \ref{res-03}]
It is an application of Theorems \ref{xP3} and \ref{T3a}. We use
the notation and framework of Section \ref{sec-saddle}.
\\
With (\ref{eq-09}) and Theorem \ref{xP3}-a,
$\dom\Fs\subset\LL=\Ll'.$ For all $\ell\in\Ll',$
\begin{eqnarray*}
  \Fos(\ell)
  &\stackrel{(a)}=&\Fs_+(\ell_+)+\Fs_-(\ell_-) \\
  &\stackrel{(b)}=& \inf\{I_{\lambda^*_+}(k^a)+\iota^*_{\dom\lambda_+}(k^s); k\in L_{\lambda_+}':k\ge0, k_{|\Ll}=\ell_+\}\\
  &&\ +\inf\{I_{\lambda^*_-}(k^a)+\iota^*_{\dom\lambda_-}(k^s); k\in L_{\lambda_-}':k\ge0, k_{|\Ll}=\ell_-\}
\end{eqnarray*}
 Equality (a) comes from Proposition
\ref{res-05}-a and equality (b) is a dual equality of the type of
Theorem  \ref{xP3}-b applied with
\begin{equation}\label{eq-23}
    I_\rho^*(k)=I_{\rho^*}(k^a)+\iota^*_{\dom\rho}(k^s)\quad k\in L_{\rho}'
\end{equation}
which holds for any Young function $\rho.$ This identity is proved
by Foug\`eres, Giner, Kozek and Rockafellar
\cite{FG76,Koz79,Roc71} under the assumptions (A$_R$) and
(A$_{\gamma^*}^1$). The function $I_\rho^*$ is strongly continuous
on $\icordom I_\rho^*\subset L_\rho',$ see \cite[Lemma
2.1]{Leo01a}. Hence, under the assumption (2), we obtain that
\begin{equation}\label{eq-20}
  \Ib(\ell)=\Fos(\ell-mR),\quad \ell\in\Ll'
\end{equation}
taking advantage of the direct sum $\ell=\oplus_k\ell_{|\ZZ_k}$
acting on $u=(u_{|\ZZ_k})_{k\ge1}$ which lead to the nonnegative
series $\Fi(u)=\oplus_k\Fi(u_{|\ZZ_k})$ and $\Fs(\ell)=\sum_k
\Fs(\ell_{|\ZZ_k}).$

 \noindent $\bullet$\ \textit{Reduction to $m=0.$}\ We have  seen at (\ref{PEnt})
that the transformation $Q\rightsquigarrow \ell=Q-mR$ corresponds
to the transformations $\gamma\rightsquigarrow\lambda$ and
(\ref{pc}) $\rightsquigarrow$ (\ref{PEnt}). This still works with
(\ref{pbc}) and one can assume from now on without loss of
generality that $m=0$ and $\gamma=\lambda.$
\\
The assumption (A$_{\gamma^*}^2$) will not be used during the rest
of the proof. This allows Remark \ref{rem-01}-d.
\\
$\bullet$\ \textit{Verification of \HF\ and \HT.}\ Suppose that
$W=\{z\in\ZZ;\lambda(z,s)=0,\forall s\in\R\}$ is such that
$R(W)>0.$ Then, any $\ell$ such that $\langle
u\mathbf{1}_W,\ell\rangle>0$ for some $u\in\Ll$ satisfies
$\Fs(\ell)=+\infty.$ Therefore, one can remove $W$ from $\ZZ$
without loss of generality. Once, this is done, the hypothesis
\HF\ is satisfied under the assumption (A$_{\gamma^*}^1$). The
hypothesis \HTi\ is (A$_\theta^\exists$) while \HTii\ is
(A$_\theta^2$).
\\
$\bullet$ \noindent\textit{The computation of $\Fb$ in the case
where $\lambda$ is even.}\
    Since $\Fi$ is even, Theorem \ref{res-01} tells us that $\dom\Fb$ is included in
    the $\sigma(L_\lambda'',L_\lambda')$-closure of $\dom\Fi.$
     Thanks to (\ref{eq-23}) and the decomposition (\ref{eq-15}), the extension $\Fb$ is given
for each $\zeta\in L_{\lambda}''$ by
\begin{eqnarray*}
  \Fb(\zeta)
  &=& (\Ib_{\lambda^*})^*(\zeta_1,\zeta_2) \\
  &=& \sup_{f\in\Lls,k\in\Ll^s}\{\langle\zeta_1,fR\rangle+\langle\zeta_2,k\rangle
  -\Ils(fR)-\iota_{\dom \Il}^*(k)\} \\
  &=& I^*_{\lambda^*}(\zeta_1)+\iota_{\dom \Il}^{**}(\zeta_2)\\
  &=& \Ib_{\lambda}(\zeta_1)+\iota_D(\zeta_2)\\
  &=& \Il(\zeta_1^a)+\iota^*_{\dom\Ils}(\zeta_1^s)+\iota_D(\zeta_2)
\end{eqnarray*}
where $D$ is the
$\sigma(L_{\lambda}^{s\prime},L_{\lambda}^s)$-closure of $\dom\Il$
and we dropped the restrictions $\zeta_{|L}$ for simplicity.
\\
$\bullet$ \noindent\textit{Extension to the case where $\lambda$
is not even.}\ By Proposition \ref{res-05}-b, we have
$\Fb(\zeta)=\Fb_+({\zeta_+}_{|L_{\lambda_+}'\cap\Ll'})+\Fb_-({\zeta_-}_{|L_{\lambda_-}'\cap\Ll'})$
if $\zeta\in\Ll'^b$ and $+\infty$ otherwise. It follows that
\begin{equation}\label{eq-12}
     \Fb(\zeta)=\Ib_\lambda^*(\zeta)
 = \Il(\zeta_1^a)+\iota^*_{\dom\Ils}(\zeta_1^s)+\iota_D(\zeta_2)
\end{equation}
if $\zeta\in K_{\lambda}''$ and $+\infty$ otherwise.
    In particular, we have
\begin{eqnarray*}
    \Lambda(y)&=&\Il(\yt),\quad y\in\YY\\
    \Lb(\omega)&=&\left\{
    \begin{array}{ll}
      \Il\big([T^*\omega]_1^a\big)+\iota^*_{\dom\Ils}\big([T^*\omega]_1^s\big)+\iota_D\big([T^*\omega]_2\big)\quad
      & \textrm{if }\omega\in\YYb \\
      +\infty & \textrm{otherwise} \\
    \end{array}
    \right.,
    \quad \omega\in\XX_L^*.
\end{eqnarray*}
This provides us with the dual problems (\ref{dc}) and
(\ref{dbc}).
\\
$\bullet$\ \textit{Proof of (a) and (b).}\ Apply Theorem
\ref{xP3}. $\Box$

Let us go on with the proof of (c).  By Theorem \ref{T3a},
(\ref{pbc},\ref{dbc}) admits a solution in $\Ll'\times\YYb$ and
$(\lh,\ob)\in\Ll'\times\YYb$ solves (\ref{pbc},\ref{dbc}) if and
only if
\begin{equation}\label{eq-37}
    \left\{\begin{array}{cl}
      (a) & \xh\in C\cap\dom\Gamma^* \\
      (b) & \langle \ob,\xh\rangle \leq \langle \ob,x\rangle, \forall x\in C\cap\dom\Gamma^* \\
      (c) & \lh\in\partial_{\Ll'}\Fb(T^*\ob) \\
    \end{array}\right.
\end{equation}
where $\xh\eqdef T\lh$ is defined in the weak sense with respect
to the duality $\langle\YY_L,\XX_L\rangle.$ Since
$\dom\Gamma^*\subset\XX_L,$ the above dual brackets are
meaningful.
\\
    $\bullet$ \noindent\textit{The computation of }$\partial_{L_\lambda'}\Fb(\zeta).$\
Let us first assume that $\lambda$ is even. For all $u\in
L_\lambda,$ $u^a_1=u_2=u$ and $u_1^s=0.$ This gives
$\Fb(\zeta+u)-\Fb(\zeta)=I_\lambda(\zeta_1^a+u_1)-I_\lambda(\zeta_1^a)+\iota_D(\zeta_2+u_2)-\iota_D(\zeta_2)$
where $u_1=u$ and $u_2=u$ act respectively on $L_{\lambda^*}R$ and
$L_\lambda^s.$ This direct sum structure leads us to
\begin{equation}\label{eq-21}
 \partial_{L_\lambda'}\Fb(\zeta)=\partial_{L_{\lambda^*}R}I_\lambda(\zeta_1^a)+\partial_{L_\lambda^s}\iota_D(\zeta_2).
\end{equation}
which again is the direct sum of the absolutely continuous and
singular components of $\partial_{L_\lambda'}\Fb(\zeta).$
Differentiating in the directions of $\UU=L_\lambda,$ one obtains
$\partial_{L_{\lambda^*}R}I_\lambda(\zeta_1^a)=\{\lambda'(\zeta_1^a)R\}.$
The computation of $\partial_{L_{\lambda}^s}\iota_D(\zeta_2)$ is
standard:
    $
\partial_{L_{\lambda}^s}\iota_D(\zeta_2)=D^\bot(\zeta_2)
    $
is the outer normal cone of $D$ at $\zeta_2.$
\\
Now, consider a general $\lambda.$ By Proposition \ref{res-05}-a,
$\lh_+\in\partial_{L_{\lambda_+}'\cap\Ll'}\Fb_+([T^*\ob]_+)$ and
$\lh_-\in\partial_{L_{\lambda_-}'\cap\Ll'}\Fb_-([T^*\ob]_-)$.
Therefore, (\ref{eq-21}) becomes
\begin{equation*}
 \partial_{\Ll'}\Fb(\zeta)=\partial_{K_{\lambda^*}R\cap\LlsR}I_\lambda(\zeta_1^a)+\partial_{K_\lambda^s\cap \Ll^s}\iota_D(\zeta_2).
\end{equation*}
$\bullet$ \noindent\textit{Representation of $[T^*\ob]_1^a.$}\ One
still has to prove that
\begin{equation}\label{eq-32}
    [T^*\ob]_1^a(z)=\langle \theta(z),\ot\rangle
\end{equation}
for $R$-a.e.\! $z\in \ZZ$ and some linear form $\ot$ on $\XXo.$
\\
 If $W_-:=\{z\in\ZZ;\lambda(z,s)=0, \forall s\le0\}$ satisfies $R(W_-)>0,$ $\dom\Ib$ is a set of
 linear forms which are nonnegative on $W$ and  $\gamma_z'(s)=0$ for all $s\le0, z\in W.$
Hence, one can take any function for the restriction to $W_-$ of
$[T^*\ob]_{1-}^a$ without modifying (\ref{eq-37})-c. As a
symmetric remark holds for $W_+=\{z\in\ZZ;\lambda(z,s)=0, \forall
s\ge0\},$ it remains to consider the situation where for
$R$-a.e.\! $z,$ there are $s_-(z)<0<s_+(z)$ such that
$\lambda(z,s_\pm(z))>0.$ This implies that
$\lim_{s\rightarrow\pm\infty}\lambda(z,s)/s>0.$
\\
By Theorem \ref{T3a}, $T^*\ob$ is in the $\sigma( K_\lambda'',
K_\lambda')$-closure of $T^\ast(\dom\La).$ Therefore,
$[T^*\ob]_1^a$ is in the $\sigma(K_\lambda,K_\lambda')$-closure of
$T^\ast(\dom\La).$ As $T^\ast(\dom\La)$ is convex, this closure is
its strong closure  in $K_\lambda.$ Since there exists a finite
measurable function $c(z)$ such that $0<c(z)\leq
\lim_{s\rightarrow\infty}\lambda(z,s)/s,$ one can consider the
nontrivial Young function $\rho(z,s)=c(z)|s|$ and the
corresponding Orlicz spaces $L_\rho$ and $L_\rho'=L_{\rho^*}.$ If
$R$ is a bounded measure, we have $L_{\lmax}\subset L_\rho$ and
$L_{\rho^*}\subset L_{\lmaxs},$ so that $[T^*\ob]_1^a$ is in the
strong closure of $T^\ast(\dom\La)$ in $L_\rho.$
\\
 As a consequence,
$[T^*\ob]_1^a$ is the pointwise limit of a sequence $(T^\ast
y_n)_{n\geq 1}$ with $y_n\in \YY.$ As $T^\ast y_n(z)=\langle
y_n,\theta(z)\rangle,$ we see that $[T^*\ob]_1^a(z)=\langle
\theta(z),\ot\rangle$ for some linear form $\ot$ on $\XXo.$ If $R$
is unbounded, it is still assumed to be $\sigma$-finite: there
exists a sequence $(\ZZ_k)$ of measurable subsets of $\ZZ$ such
that $\cup_k \ZZ_k=\ZZ$ and $R(\ZZ_k)<\infty$ for each $k.$ Hence,
for each $k$ and all $z\in \ZZ_k,$ $(T^*\ob)^a(z)=\langle
\theta(z),\ot^k\rangle$ for some linear form $\ot^k$ on $\XXo,$
from which (\ref{eq-32}) follows.
\\
$\bullet$\noindent\textit{Proof of (c).} It follows from the
previous considerations and Theorem \ref{T3a}.
\\
$\bullet$\noindent\textit{Proof of (d).}\ Statement (d)-1 follows
from Theorem \ref{xP3}. Statement (d)-2 is immediately deduced
from (c). Finally, (d)-3 is (\ref{xeq-11}).
\end{proof}

\section{Solving ($\pc$)}\label{sec-pc}

The general assumptions (A) are imposed and we study (\ref{pc})
under the additional good constraint assumption (\ref{A-forall})
which imposes that the convex set $C$ is such that
 \begin{equation}\label{eq-05}
   \To^{-1}C\cap\LlsR=\bigcap_{y\in Y}\left\{fR\in\LlsR;\IZ \yt f\,dR\ge a_y\right\}
\end{equation}
for some subset $Y\in\XXo^*$ such that $\yt\in \El$ for all $y\in
Y$ and some function $y\in Y\mapsto a_y\in\R.$
\\
The dual problem ($\mathcal{D}$) associated with (\ref{pc}) is
(\ref{dc}) and the extended dual problem is
\begin{equation}\label{dtc}
    \textsl{maximize } \inf_{x\in C}\langle\omega,x\rangle-I_\gamma(\langle\omega,\theta\rangle),
   \quad \omega\in\YYt\tag{$\dtc$}
\end{equation}
where $\YYt$ is the convex cone of all linear forms $\omega$ on
$\XXo$ which are such that
\begin{enumerate}[-]
    \item the function $\langle\omega,\theta(\cdot)\rangle_{\XXo^*,\XXo}$ is
    measurable;
    \item
    $\IZ\lambda(t\langle\omega,\theta(\cdot)\rangle)\,dR<\infty$
    for some $t>0;$
    \item $\langle\omega,\theta(\cdot)\rangle$ is in the
    $\sigma(K_\lambda,K_{\lambda^*})$-closure of $\{\yt;y\in\YYo\}.$
\end{enumerate}

\begin{theorem}\label{res-02}
Suppose that
\begin{enumerate}
    \item the assumptions (A) and (\ref{A-forall}) are satisfied;
    \item for $R$-almost every $z\in \ZZ,$
    $ \lim_{t\rightarrow \pm\infty}\gamma_z^*(t)/t=+\infty;$ % eq-07
    \item $C$ satisfies (\ref{eq-05}) with $\yt\in\El$ for all $y\in Y.$
\end{enumerate}
Then:
\begin{enumerate}[(a)]
 \item The dual equality for \emph{(\ref{pc})} is
    \begin{equation*}
    \inf(\pc)=\sup(\dc)=\sup(\dtc)=\inf_{x\in C}\Gamma^*(x)\in [0,\infty].
    \end{equation*}
 \item If $C\cap \dom\Gamma^*\not=\emptyset$ or equivalently $C\cap \To\dom I\not=\emptyset,$
 then \emph{(\ref{pc})}
 admits a unique solution $\Qh$ in $\LlsR$ and any minimizing sequence $\seq Qn$ converges to
 $\Qh$ with respect to the topology   $\sigma(L_{\lmaxs}.R,\El).$
\end{enumerate}
 Suppose that in addition $C\cap \icordom\Gamma^*\not=\emptyset$ or equivalently $C\cap \icor(\To\dom I)\not=\emptyset.$
\begin{enumerate}[(a)]
 \item[(c)] Let us define $\xh\eqdef \IZ\theta\,d\Qh$ in the weak sense with respect to
    the duality $\langle\YYo,\XXo\rangle.$
 There exists  $\ot\in\YYt$ such that
 \begin{equation}\label{eq-36}
    \left\{\begin{array}{cl}
      (a) & \xh\in C\cap\dom\Gamma^* \\
      (b) & \langle \ot,\xh\rangle_{\XXo^*,\XXo} \leq \langle \ot,x\rangle_{\XXo^*,\XXo}, \forall x\in C\cap\dom\Gamma^* \\
      (c) & \Qh(dz)=\gamma'_z(\langle\ot, \theta(z)\rangle)\,R(dz). \\
    \end{array}\right.
\end{equation}
Furthermore, $\Qh\in\LlsR$ and $\ot\in\YYt$ satisfy (\ref{eq-36})
if and only if $\Qh$ solves \emph{(\ref{pc})} and $\ot$ solves
\emph{(\ref{dtc})}.
    \item[(d)]
    Of course, (\ref{eq-36}-c) implies
    $
\xh=\IZ \theta\gamma'(\ott)\,dR
    $
in the weak sense.
    Moreover,
    \begin{enumerate}[1.]
        \item $\xh$ minimizes $\Gamma^*$ on $C,$
        \item
        $I(\Qh)=\Gamma^*(\xh)=\IZ\gamma^*\circ\gamma'(\ott)\,dR<\infty$
        and
        \item $I(\Qh)+\IZ\gamma(\ott)\,dR=\IZ\ott\,d\Qh.$
    \end{enumerate}
\end{enumerate}
\end{theorem}

\begin{proof}
It is a corollary of the proof of Theorem \ref{res-03}.  One
applies the abstract results of Section \ref{sec-saddle} with
\begin{equation}\label{eq-18}
     \Fo(u)=\Il(u):= \IZ\lambda(u)\,dR,\quad u\in\UUo:=\EEl
\end{equation}
This gives $\UU=\El$ with the Orlicz norm $|u|_\Fi=\|u\|_\lmax$
and $\LL=\LlsR.$ The space $\YY:=\YY_E$  is the completion of
$\YYo$ endowed with the norm $|y|_\Lambda=\|\yt\|_\lmax.$ It is
isomorphic to the closure of the subspace $\{\yt;y\in\YYo\}$ in
$\El,$ see assumption (\ref{A-forall}). The topological dual space
$\XX_E=\YY_E'$ is identified with $\LlsR/\mathrm{ker\,}T$ and its
norm is given by $|x|_\Lambda^*=\inf\{\|f\|_\lmaxs;f\in\Lls:
T(fR)=x\}.$
\\
The assumption (3) is: $C$ is a convex
$\sigma(\XX_E,\YY_E)$-closed set.
\\
As in the proof of Theorem \ref{res-03}, one reduces to the case
where $m=0$ without loss of generality.
\\
The assumption (2) implies that $\lambda$ is a finite function. It
follows that $\El'=\Lls,$ the convex conjugate $\Fs$ of $\Fi$ with
respect to the duality $\langle\El,\Lls\rangle$ is
$$\Fs=I_{\lambda^*}$$ (see \cite{Roc68}) and the corresponding extended function
$\Fb$ is
    \begin{equation*}
  \Fb(\zeta)=\Il(\zeta^{a})+\iota_{\dom \Ils}^*(\zeta^s)
    \end{equation*}
    if $\zeta$ is in $K_{\lambda}.R\oplus K_{\lambda}^s$ and $+\infty$
    otherwise.
\\
With these correspondences, the proof of the theorem is an
immediate translation of the proof of Theorem \ref{res-03}.
\end{proof}

\begin{remarks} \
\begin{enumerate}[(a)]
    \item
    The assumption (2) implies that $\lambda$ is a
    finite function. Note that otherwise one would get $\El=\{0\}.$
    \item
    As in Remark \ref{rem-01}-d, removing the assumption
(A$_{\gamma^*}^2$): $m\in\Lls,$ one can still consider the
minimization problem (\ref{pco}) instead of (\ref{pc}). The
transcription of Theorem \ref{res-02} is as follows. Replace
respectively (\ref{pc}), $C,$ $\Gamma^*,$ $\xh$ and $\gamma$ by
(\ref{pco}), $\Co,$ $\Ls,$ $\tilde{x}$ and $\lambda$ where
$\tilde{x}=\IZ\theta\,d(\Qh-mR)$ is well-defined.
\\
The statement (b) must be replaced by the following one:
\textit{If $\Co\cap \dom\Lambda^*\not=\emptyset,$ then
\emph{(\ref{pco})}
 admits a unique solution $\Qh$ in $mR+\LlsR$ and any minimizing sequence $\seq
    Qn$ is such that $(Q_n-mR)_{n\ge1}$
    converges in $\LlsR$ to $\Qh-mR$ with respect to the topology $\sLE.$}
    \item Seeing Theorem \ref{res-02} as a direct corollary of
    Theorem \ref{res-03} would have been possible since
    Proposition \ref{res-B4} insures that $\To^{-1}C\cap\Ll'
    =\bigcap_{y\in Y}\left\{\ell\in\Ll';\langle\theta,\ell\rangle\ge a_y\right\}
    =\bigcap_{y\in Y}\left\{fR\in\LlsR;\IZ \yt f\,dR\ge
    a_y\right\}$ whenever (\ref{A-forall}) holds. But, the
    drawback is that the unnecessary assumption (2) of Theorem
    \ref{res-03} has to be kept.
\end{enumerate}
\end{remarks}

\section{Examples}

Standard examples of entropy minimization problems are presented.

\subsection{Some examples of entropies}\label{sec-expl-ent}
The entropies defined below occur naturally in statistical
physics, probability theory, mathematical statistics and
information theory.

\subsubsection*{Boltzmann entropy} The Boltzmann entropy with respect to the positive measure
$R$ is defined by
    $
    H_B(Q|R)=\left\{%
\begin{array}{ll}
    \IZ \log\left(\frac{dQ}{dR}\right)\,dQ & \hbox{if }0\le Q\prec R \\
    +\infty, & \hbox{otherwise} \\
\end{array}%
\right.
    $
for each $Q\in\MZ.$
It corresponds to $\gamma^*_z(t)=\left\{%
\begin{array}{ll}
    t\log t & \hbox{if }t>0 \\
     0 & \hbox{if }t=0 \\
    +\infty & \hbox{if }t<0 \\
\end{array}%
\right..$ But this $\gamma^*$ takes negative values and is ruled
out by our assumptions. A way to circumvent this problem is to
consider the variant below.

\subsubsection*{A variant of the Boltzmann entropy} Let $m:\ZZ\to(0,\infty)$
be a positive measurable function. Considering
$$\gamma^*_z(t)=t\log t -[1+\log m(z)]t+m(z),\quad t>0,$$ one sees that
it is nonnegative and that $\gamma^*_z(t)=0$ if and only if
$t=m(z).$ Hence $\gamma^*$ enters the framework of this paper and
\begin{equation}\label{eq-02}
    \lambda_z(s)=m(z)[e^s-s-1],\quad s\in\R.
\end{equation}
It is easily seen that
\begin{equation*}
    H_B(Q|R)=I_{\gamma^*}(Q)+\IZ (1+\log m)\,dQ-\IZ m\,dR
\end{equation*}
which is meaningful if $Q$ integrates $1+\log m$ where  $m\in
L^1(R).$
\\
As an application, let $R$ be the Lebesgue measure on $\ZZ=\R^d$
and minimize $H_B(Q|R)$  on the set $\CC=\{Q\in\PZ; \IZ
|z|^2\,Q(dz)=E\}\cap \CC_o.$ Taking $m(z)=e^{-|z|^2},$ one is led
to minimizing $I_{\gamma^*}$ on $\CC.$

\subsubsection*{A special case} It is defined by
\begin{equation}\label{eq-28}
    H(Q|R)=\left\{%
\begin{array}{ll}
    \IZ \left[\frac{dQ}{dR}\log\left(\frac{dQ}{dR}\right)-\frac{dQ}{dR}+1\right]\,dR & \hbox{if }0\le Q\prec R \\
    +\infty, & \hbox{otherwise} \\
\end{array}%
\right.,\quad Q\in\MZ
\end{equation}
It corresponds to $\gamma^*_z(t)=\left\{%
\begin{array}{ll}
    t\log t-t+1 & \hbox{if }t>0 \\
     1 & \hbox{if }t=0 \\
    +\infty & \hbox{if }t<0 \\
\end{array}%
\right.,   $ $m(z)=1$ and
    $
    \lambda_z(s)=e^s-s-1,$ $s\in\R
    $
 for all $z\in \ZZ.$  Note that
$H(Q|R)<\infty$ implies that $Q$ is nonnegative.

\subsubsection*{Relative entropy}  The reference measure $R$ is assumed to be
a probability measure and one denotes $\PZ$ the set of all
probability measures on $\ZZ.$ The relative entropy of $Q\in\MZ$
with respect to $R\in\PZ$ is the following variant of the
Boltzmann entropy:
\begin{equation*}
    I(Q|R)=\left\{%
\begin{array}{ll}
    \IZ \log\left(\frac{dQ}{dR}\right)\,dQ & \hbox{if }Q\prec R \hbox{ and }Q\in\PZ\\
    +\infty & \hbox{otherwise} \\
\end{array}%
\right.,\quad Q\in\MZ.
\end{equation*}
It is (\ref{eq-28}) with the additional constraint that
$Q(\ZZ)=1:$
$$
I(Q|R)=H(Q|R)+\iota_{\{Q(\ZZ)=1\}}
$$
When minimizing the Boltzmann entropy $Q\mapsto H_B(Q|R)$ on a
constraint set which is included in $\PZ,$ we have for all
$P,Q\in\PZ,$
\begin{equation*}
    H_B(Q|R)= I(Q|P)+\IZ \log\left(\frac{dP}{dR}\right)\,dQ
\end{equation*}
which is meaningful for each $Q\in\PZ$ which integrates
$\frac{dP}{dR}.$

\subsubsection*{Reverse relative entropy} The reference measure $R$ is assumed to be
a probability measure. The reverse relative entropy is
\begin{equation*}
    Q\in\MZ\mapsto \left\{%
\begin{array}{ll}
    I(R|Q) & \hbox{if }Q\in\PZ \\
    +\infty & \hbox{otherwise}  \\
\end{array}%
\right.     \in [0,\infty].
\end{equation*}
It corresponds to
$\gamma^*_z(t)=\left\{%
\begin{array}{ll}
    -\log t +t-1 & \hbox{if }t>0 \\
    +\infty & \hbox{if }t\le 0 \\
\end{array}%
\right.,$
 $m(z)=1$   and
\begin{equation}\label{eq-03}
    \lambda_z(s)=
\left\{%
\begin{array}{ll}
    -\log(1-s)-s & \hbox{if }s<1 \\
    +\infty & \hbox{if }s\ge1 \\
\end{array}%
\right.,
\end{equation}
for all $z\in\ZZ,$ with the additional constraint that $Q(\ZZ)=1.$

\subsection{Some examples of constraints}\label{sec:constraints}

Let us consider two standard constraints which are the moment
constraints and the marginal constraints.

\subsubsection*{Moment constraints}

\newcommand{\AB}{{A\!\times\! B}}
\newcommand{\IAB}{\int_{\AB}}

Let $ \theta=(\theta_k)_{1\leq k\leq K} $ be a measurable function
from $\ZZ$ to  $\XXo=\R^K.$ The moment constraint is specified by
the operator
$$
\To\ell=\IZ\theta\,d\ell=\left(\IZ\theta_k\,d\ell\right)_{1\leq
k\leq K}\in\R^K,
$$
which is defined for each $\ell\in\MZ$ which integrates all the
real valued measurable functions $\theta_k.$ The adjoint operator
is
\begin{equation*}
    \To^\ast y(z)=\sum_{1\leq k\leq K} y_k\theta_k(z),\quad
    y=(y_1,\dots,y_K)\in\R^K, z\in\ZZ.
\end{equation*}

\subsubsection*{Marginal constraints}
Let $\ZZ=\AB$ be a product space, $M_{AB}$ be the space of all
\emph{bounded} signed measures on $\AB$ and $U_{AB}$ be the space
of all measurable bounded functions $u$ on $\AB.$ Denote
$\ell_A=\ell(\cdot\times B)$ and $\ell_B=\ell(A\times\cdot)$ the
marginal measures of $\ell\in M_{AB}.$ The constraint of
prescribed marginal measures is specified by
$$
\IAB\theta\,d\ell=(\ell_A,\ell_B)\in M_{A}\times M_{B},\quad
\ell\in M_{AB}
$$
where $M_{A}$ and $M_{B}$ are the spaces of all bounded signed
measures on $A$ and $B.$ The function $\theta$ which gives the
marginal constraint is
$$
\theta(a,b)=(\delta_{a}, \delta_{b}),\ a\in A, b\in B
$$
where $\delta_a$ is the Dirac measure at $a.$ Indeed,
$(\ell_A,\ell_B)=\IAB (\delta_a,\delta_b)\,\ell(dadb).$
\\
More precisely, let $U_{A},$ $U_{B}$ be the spaces of measurable
functions on $A$ and $B$ and take $\YYo=U_{A}\times U_{B}$ and
$\XXo=U_{A}^*\times U_{B}^*.$ Then, $\theta$ is a measurable
function from $\ZZ=\AB$ to $\XXo=U_{A}^*\times U_{B}^*.$ It is
easy to see that the adjoint of the marginal operator
\begin{equation*}
    \To\ell=(\ell_A,\ell_B)\in U_{A}^*\times U_{B}^*,\quad
    \ell\in\LLo=U_{AB}^*
\end{equation*}
where $\langle f,\ell_A\rangle:=\langle f\otimes 1,\ell\rangle$
and $\langle g,\ell_B\rangle:=\langle 1\otimes g,\ell\rangle$ for
all $f\in U_A$ and $g\in U_B,$ is given by
\begin{equation}\label{eq-19}
    \To^\ast(f,g)=f\oplus g\in U_{AB},\quad f\in U_{A}, g\in U_{B}
\end{equation}
where $f\oplus g(a,b):=f(a)+g(b),$ $a\in A, b\in B.$

%\bibliographystyle{plain}
%\bibliography{bib-christian}

\end{document}